\documentclass[a4paper,draft]{amsart}
\usepackage{amsmath,amssymb,amscd,amsthm}

\newtheorem{formula}{}[section]

\newtheorem{corollary}[formula]{Corollary}

\newtheorem{theorem}[formula]{Theorem}
\theoremstyle{definition}
\newtheorem{definition}[formula]{Definition}

\theoremstyle{remark}
\newtheorem*{remark}{Remark}

\renewcommand{\a}{\alpha}
\newcommand{\g}{\gamma}
\renewcommand{\l}{\lambda}
\renewcommand{\L}{\Lambda}
\newcommand{\Z}{\mathbb Z}
\newcommand{\C}{\mathbb C}
\newcommand{\Q}{\mathbb Q}
\newcommand{\m}{\mu}
\newcommand{\n}{\nu}
\newcommand{\s}{\sigma}
\renewcommand{\O}{\Omega}
\renewcommand{\o}{\omega}

\begin{document}

\title[Cobordism classification of manifolds with $\Z/p$-action]
{On the cobordism classification of manifolds with $\Z/p$-action}
\author{Taras E. Panov}
\thanks{Partially
supported by Russian Foundation of Fundamental Research
grant no. 99-01-00090.}
\address{Department of Mathematics and Mechanics, Moscow State
University, 119899 Moscow, Russia}
\email{tpanov@mech.math.msu.su}

\begin{abstract}
We refer to an action of the group $\Z/p$ ($p$ here is an odd prime) on a
stably complex manifold as {\it simple} if all its fixed
submanifolds have the trivial normal bundle. The important particular case of
a simple action is an action with only isolated fixed points. The problem of
cobordism classification of manifolds with simple action of $\Z/p$ was posed
by V.\,M.~Buchstaber and S.\,P.~Novikov in 1971. The analogous question of
cobordism classification with stricter conditions on $\Z/p$-action was
answered by Conner and Floyd. Namely, Conner and Floyd solved the problem in
the case of simple actions with identical sets of weights (eigenvalues of the
differential of the map corresponding to the generator of $\Z/p$) for all
fixed submanifolds of same dimension. However, the general setting of the
problem remained unsolved and is the subject of our present paper. We have
obtained the description of the set of cobordism classes of stably complex
manifolds with simple $\Z/p$-action both in terms of the coefficients of
universal formal group law and in terms of the characteristic numbers, which
gives the complete solution to the above problem. In particular, this gives
a purely cohomological obstruction to the existence of a simple
$\Z/p$-action (or an action with isolated fixed points) on a manifold. We also
review connections with the Conner--Floyd results and with the well-known
Stong--Hattori theorem.
\end{abstract}

\maketitle

\section{Introduction}

By {\it a stably (or weakly almost) complex manifold} we mean an orientable
differentiable manifold with complex structure in its stable tangent
bundle~--- the main object of study in the complex cobordism theory. The
fixed point set of a smooth action of $\Z/p$ on a differentiable manifold
$M$ (i.e. the fixed point set of a diffeomorphism of an odd prime period) is
decomposed into the sum of connected fixed submanifolds of even
codimensions~(cf.~\cite{CF1}). The action of $\Z/p$ induces a representation
of $\Z/p$ in the tangent space $T_qM$ at each fixed point $q\in M$. This
representation is decomposed into the sum of non-trivial one-dimensional
complex representations of $\Z/p$ and the trivial representation in the
tangent subspace to the fixed submanifold containing $q$. This decomposition
is same for all fixed points within the same connected fixed submanifold.
Each non-trivial one-dimensional representation of $\Z/p$ is determined by an
element $x_k\in(\Z/p)^{*}$, $x_k\ne0\mod p$ (the generator of $\Z/p$
acts by multiplication by $e^{2\pi ix_k/p}$). We refer to $x_1,\ldots,x_m$
as {\it the weights} of the $\Z/p$-action on $M$ corresponding to the fixed
submanifold (here $m$ equals half the codimension of the fixed submanifold in
$M$).

\vspace{-2pt}

\begin{definition}
\label{simple}
  We refer to an action of $\Z/p$ on a stably complex
  manifold $M$ as {\it simple} if all the
  fixed submanifolds have trivial normal bundle.
  A simple action is called {\it strictly simple},
  if the sets of weights for the action are identical for all fixed
  submanifolds of the same dimension.
\end{definition}

\vspace{-2pt}

The good example of a simple action is an action with only isolated
fixed points.

In the present article we obtain the full classification of complex
cobordism classes $\s\in\O_U$ containing a manifold with simple action of
$\Z/p$. The description is given both in terms of the coefficients of
universal formal group law of ``geometric cobordisms" (theorem~\ref{form})
and in terms of the characteristic numbers (theorem~\ref{charnum},
corollary~\ref{cobord}). Our corollary~\ref{cobord} can be also
viewed as a cohomological obstruction to the existence of a simple
$\Z/p$-action (or an action with isolated fixed points) on a manifold.

The classification problem for strictly simple actions of
$\Z/p$ was completely solved by Conner and Floyd in~\cite{CF1}.
(A strictly simple action from definition~\ref{simple} was called
in~\cite{CF1} just ``an action of $\Z/p$ with fixed point set having
the trivial normal bundle".) Note, that even in the special case of an action
with finite number of isolated fixed points the notions of simple and
strictly simple action differ (see examples below).  The Conner--Floyd
results follow from the results of our paper. At the same time, it seems to
us that the approach used in~\cite{CF1} does not allow to obtain our more
general result.

The applications of the formal group law theory to
$\Z/p$-actions were firstly discussed in the pioneer article~\cite{No1}.
The formal group law theory itself arises in topology due to so-called
formal group law of geometric cobordisms. The problem solved in our article
was firstly formulated in~\cite{BuN}. There was obtained the formula
expressing the $\mod p$ cobordism class of a manifold $M$ with a simple
action of $\Z/p$ in terms of the cobordism classes of fixed submanifolds and
the weights of the action (see formula~(\ref{class})). Actually, the first
results on the problem were obtained even earlier, in~\cite{Ka}.
Particularly, there was proved the statement mentioned in our article as
corollary~\ref{smdim}. In \cite{Ka}, as well as in our work, the set of
cobordism classes of manifolds with a simple $\Z/p$-action is handled as
the $\O_U$-module spanned by certain coefficients of the power system
determined by the formal group law of geometric cobordisms.  (Here $\O_U$ is
the complex cobordism ring of point, which is isomorphic to the polynomial
ring $\Z[a_1,a_2,\ldots]$, $\deg a_i=-2i$, as it was shown by Milnor and
Novikov.) However, our new choice of generators for the above $\O_U$-module
allowed us to present the answer for the classification problem in terms of
the characteristic numbers, which of course is much more explicit.

\section{The ring $U^{*}(\Z/p)$ of equivariant cobordisms
with free $\Z/p$-action and the Conner--Floyd equations}

Now, let us start with a stably complex manifold $M$ with a simple
action of $\Z/p$. Let the fixed submanifolds represent the
cobordism classes $\l_j\in\O_U$ and have weights
$(x_1^{(j)},\ldots,x_{m_j}^{(j)})$ in their (trivial) normal
bundles. Here $2m_j+\dim\l_j=\dim M$,
$(x_k^{(j)})\in(\Z/p)^*$, $k=1,\ldots,m_j$. These data define the
cobordism class of $M$ in $\O_U$ up to elements from $p\O_U$
(cf.~\cite{BuN}). This follows from the fact that the cobordism class of
a manifold with free $\Z/p$-action (i.e. without fixed points) belongs to
$p\O_U$, and vice versa, any cobordism class from $p\O_U$ can be represented
by a manifold with free action of $\Z/p$.

Let $B\Z/p$ denote the classifying space for the group $\Z/p$
(the infinite-dimensional lens space).
To each fixed submanifold of the cobordism class $\l_j\in\O_U$ with weights
$(x_1^{(j)},\ldots,x_{m_j}^{(j)})$, $2m_j+\dim\l_j=\dim M$, one can
associate an element
$$
  \alpha_{2m_j-1}(x_1^{(j)},\ldots,x_{m_j}^{(j)})\in
  U_{2m_j-1}(B{\mathbb Z}/p)
$$
called the Conner--Floyd invariant (cf.~\cite{No1}). To define
it we mention that the bordism group $U_{*}(B\Z/p)$ is isomorphic to the
group of free $\Z/p$-equivariant bordisms (cf.~\cite{CF1}). This
isomorphism takes an equivariant bordism class represented by a manifold $N$
with free action of $\Z/p$ to the bordism class in $U_{*}(B\Z/p)$ determined
by the classifying map $N/(\Z/p)\to B\Z/p$. Then
$\alpha_{2m_j-1}(x_1^{(j)},\ldots,x_{m_j}^{(j)})$ is defined as the
equivariant bordism class of the unit sphere in the fibre of (trivial) normal
bundle to $\l_j$. If we write this unit sphere as
$$
\bigl\{(z_1,\ldots,z_{m_j})\in\C^{m_j}:
|z_1|^2+\ldots+|z_{m_j}|^2=1\bigr\},
$$
then the generator of $\Z/p$ acts on it as
$$
  (z_1,\ldots,z_{m_j})\to
  \bigl(e^{2\pi ix_1^{(j)}/p}z_1,\ldots,e^{2\pi ix_{m_j}^{(j)}/p}z_{m_j}
  \bigr).
$$
The complex cobordism ring of $\C P^{\infty}$ is $\O_U[[v]]$, where
$v=c_1^U(\zeta)$ is the first cobordism Chern class of the universal line
bundle $\zeta$ over $\C P^{\infty}$ (this follows in the standard way from
the Atiyah--Hirzebruch spectral sequence in cobordisms). There is a canonical
bundle $f:B\Z/p\to\C P^{\infty}$ with fibre $S^1$. The complex cobordism ring
of $B\Z/p$ is $U^{*}(B\Z/p)=\O_U[[u]]/([u]_p=0)$ (cf.~\cite{BuN}), where
$u=f^{*}(v)$ and $[u]_p=pu+\ldots$ is the $p$-th power in the (universal)
formal group law of geometric cobordisms (a series with coefficients in
$\O_U$). We have
$D\bigl(\a_i(1,\ldots,1)\bigr)=u^{n-i}$, where $D$ is the Poincar\'e--Atiyah
duality operator from $U_*(L^{2n-1}_p)$ to $U^*(L^{2n-1}_p)$
(here $L^{2n-1}_p$ is the standard $(2n-1)$-dimensional lens space).
It follows that the $\O_U$-module $\widetilde{U}_*(B\Z/p)$
is generated by the elements $\a_{2i-1}(1,\ldots,1)$
with the following relations:
\begin{equation}
\label{rel1}
  0=\frac{[u]_p}u\cap\a_{2i-1}(1,\ldots,1).
\end{equation}
Here $\cap$ denotes the cobordism $\cap$-product, and
$$
  u^k\cap\a_{2i-1}(1,\ldots,1)=\a_{2(i-k)-1}(1,\ldots,1).
$$

It was shown in \cite{No2} that
\begin{equation}
\label{rel2}
  \alpha_{2k-1}(x_1,\ldots,x_k)=\Bigl( \prod_{j=1}^k\frac
  u{[u]_{x_j}} \Bigr) \cap \alpha_{2k-1}(1,\ldots,1).
\end{equation}
Here and below $[u]_k$ denotes the $k$-th power in the formal group law of
geometric cobordisms.
Adding the elements $\a_{2k-1}(x_1,\ldots,x_k)$, $x_j\ne 1\mod p$, to the
set of generators for the module $\widetilde{U}_*(B\Z/p)$ and adding
relations~(\ref{rel2}) to relations~(\ref{rel1}) we obtain the
$\O_U\otimes\Z_{(p)}$-free resolution of the module
$\widetilde{U}_*(B\Z/p)$ (here $\Z_{(p)}$ is a ring of
rational numbers whose denominators are relatively prime with $p$, i.e. the
ring of integer $p$-adics):
$$
  0\longrightarrow F_1\longrightarrow F_0\longrightarrow
  \widetilde{U}_*(B\Z/p)\longrightarrow 0.
$$
Here $F_0$ is the free $\O_U\otimes\Z_{(p)}$-module spanned
by $\a_{2k-1}(x_1,\ldots,x_k)$, and $F_1$ is the free
$\O_U\otimes\Z_{(p)}$-module spanned by the elements
\begin{align*}
  a(x_1,\ldots,x_k) &= \a_{2k-1}(x_1,\ldots,x_k)-\Bigl( \prod_{j=1}^k\frac
  u{[u]_{x_j}} \Bigr) \cap \alpha_{2k-1}(1,\ldots,1),\\
  a_k &= \frac{[u]_p}u\cap\a_{2k-1}(1,\ldots,1),
\end{align*}
(see formulae~(\ref{rel1}) and~(\ref{rel2})).

Hence, a simple action of $\Z/p$ on $M$ gives rise to a certain relation
between the elements $\a_{2k-1}(x_1,\ldots,x_k)$ in
$\widetilde{U}_*(B\Z/p)$. Since any element from
$\widetilde{U}_*(B\Z/p)$ corresponds to a bordism class of manifold
with free $\Z/p$-action, the converse is also true: any relation in
$\widetilde{U}_*(B{\Z}/p)$ of the form
$$
  \sum_j\l_j\a_{2m_j-1}(x_1^{(j)},\ldots x_{m_j}^{(j)})=0,\quad
  \l_j\in\O_U,\quad 2m_j+\dim\l_j=\dim M,
$$
is realized on a certain manifold $M$ with a simple
action of $\Z/p$ whose cobordism class in $\O_U$ is uniquely
determined up to elements from $p\O_U$. This manifold $M$ can be
constructed as follows. The relation in $\widetilde{U}_*(B\Z/p)$ gives us
a manifold with free $\Z/p$-action whose boundary is the union of manifolds
of the form $\l_j\times S^{2m_j-1}$. Then we glue ``covers" of the form
$\l_j\times D^{2m_j}$ to the boundary to get the closed manifold
$M$, which realizes the above relation.
Therefore, one can define the ``realization homomorphism"
\begin{equation}
\label{realiz}
  \Phi:F_1\to\O_U/p\O_U=\O_U\otimes\Z/p.
\end{equation}
It takes a relation between the elements
$\a_{2k-1}(x_1,\ldots,x_k)\in\widetilde{U}_*(B\Z/p)$ to the $\mod p$
cobordism class of the manifold that realizes this relation as
described above. The Conner--Floyd results
(cf.~\cite{CF1}, see also~\cite{Ka}) give us the following values of $\Phi$
on the basis relations from $F_1$:
\begin{align*}
  \Phi\bigl(a(x_1,\ldots,x_k)\bigr) &=
  \Bigl\langle \prod_{i=1}^k\frac u{[u]_{x_i}}
  \Bigr\rangle_k\!\mod p\:\in\:\O_U/p\O_U,\\
  \Phi(a_k)&=
  -\Bigl\langle \frac {[u]_p}u \Bigr\rangle_k\!\mod p\:\in\:\O_U/p\O_U,
\end{align*}
where $\langle\;\rangle_k$ stands for the coefficient of $u^k$.
Following \cite{BuN}, let us consider the coefficient ring $\L(1)\subset\O_U$
of the power system $[u]_k$ (i.e. the subring of $\O_U$ generated by the
coefficients of the series $[u]_k$ for all $k>0$), and its positive part
$\L^{+}(1)$, consisting of elements of non-zero degree in $\O_U$.
It follows that $\mathop{{\rm Im}}\Phi=\widetilde{\L}(1)\otimes\Z/p$, where
$\widetilde{\L}(1)=\L^+(1)\cdot\O_U$ is the $\O_U$-module spanned by
$\L^+(1)$. The homomorphism $\Phi$ lifts to
the homomorphism $F_1\to\widetilde{\L}(1)\otimes\Z_{(p)}$, or to the
homomorphism $F_1\to\widetilde{\L}_p(1)\otimes\Z_{(p)}$, where
$\widetilde{\L}_p(1)$ is the $\O_U$-module spanned by $\L^+(1)$ and $p$.
Both these homomorphisms will be denoted by the same letter $\Phi$.

Thus, the problem of description of cobordism classes of manifolds
with a simple $\Z/p$-action is equivalent to the problem of description
of $\O_U\otimes\Z/p$-module $\widetilde{\L}(1)\otimes\Z/p$, or
$\O_U\otimes\Z_{(p)}$-modules $\widetilde{\L}(1)\otimes\Z_{(p)}$,
$\widetilde{\L}_p(1)\otimes\Z_{(p)}$. These modules are
ideals in $\O_U\otimes\Z/p$ and $\O_U\otimes\Z_{(p)}$ correspondingly.

\section{The generator sets for $\O_U\otimes{\mathbb Z}/p$-module
$\widetilde{\L}(1)\otimes{\mathbb Z}/p$ and
$\O_U\otimes{\mathbb Z}_{p}$-modules
$\widetilde{\L}(1)\otimes{\mathbb Z}_{p}$,
$\widetilde{\L}_p(1)\otimes{\mathbb Z}_{p}$.}

Let us write
$$
  [u]_k=ku+\sum_{n\ge 1}\a_n^{(k)}u^{n+1},
$$
i.e. $\a_n^{(k)}\in\O_U^{-2n}$ are the coefficients of the power system.
The module $\widetilde{\L}(1)\otimes\Z_{(p)}$ is therefore generated by
$\a_n^{(k)}$, $k>1$, $n\ge1$, over $\O_U\otimes\Z_{(p)}$. The following
theorem shows that the generator set can be taken in such a way that there is
only one generator in each dimension.
\begin{theorem}
\label{form}
  One can take the following coefficients $\a_n\in\O_U^{-2n}$
  as generators of the $\O_U\otimes\Z_{(p)}$-module
  $\widetilde{\L}(1)\otimes\Z_{(p)}$:
  \begin{equation}
  \label{3}
  \a_n=\left\{
  \begin{aligned}
    \a_n^{(p_1)}, &\quad \text{if $n$ is not divisible by $p-1$},\\
    \a_n^{(p)}, &\quad \text{if $n$ is divisible by $p-1$},
  \end{aligned}
  \right.\quad n=1,2,\ldots,
  \end{equation}
  Here $p_1$ is any prime generator of the cyclic group $(\Z/p)^*$.
\end{theorem}

\begin{remark}
It follows from the Dirichlet theorem that one can
choose a {\it prime} generator $p_1$ of the cyclic group $(\Z/p)^*$.
\end{remark}

\begin{proof}[Proof of theorem~\ref{form}.]
First, let us consider
the coefficients $\a_n^{(r)}$ for non-prime $r$.
So, let $r=p_1q$ with prime $p_1$. Since
$[x]_r=[[x]_{p_1}]_q$, we have:
\begin{align*}
  rx+\sum_n\a_n^{(r)}x^{n+1}
  &=q[x]_{p_1}+\sum_n\a_n^{(q)}([x]_{p_1})^{n+1}\\
  &=p_1qx+q\sum_n\a_n^{(p_1)}x^{n+1}+\sum_n\a_n^{(q)}
  \Bigl(p_1x+\sum_m\a_m^{(p_1)}x^{m+1}\Bigr)^{n+1}.
\end{align*}
Taking the coefficient of
$x^{r+1}$ in both sides of the above identity, we get
$$
  \a_n^{(r)}=P(\a_1^{(p_1)},\ldots,\a_n^{(p_1)},\a_1^{(q)},\ldots,
  \a_n^{(q)}),
$$
where $P$ is a polynomial with integer coefficients (and zero
constant term). Hence, we can write
$$
  \a_n^{(r)}=\l_1\a_1^{(p_1)}+\ldots+\l_n\a_n^{(p_1)}+\m_1\a_1^{(q)}+
  \ldots+\m_n\a_n^{(q)},\quad\l_i,\m_i\in\O_U.
$$
Therefore, the coefficients $\a_n^{(r)}$, $r=p_1q$ can be excluded
from the set of generators for the
$\O_U\otimes\Z_{(p)}$-module $\widetilde{\L}(1)\otimes\Z_{(p)}$. Now, if
$q$ is still not prime, we repeat the above procedure until we arrive
at a set of generators consisting only of
coefficients $\a_n^{(p_1)}$ with prime $p_1$.
Now, what we need to show is that this
set of generators can still be reduced to the set~(\ref{3}).

Note, that for any (prime) generator $p_1$ of the cyclic group $(\Z/p)^*$ one
can take the coefficient $\a_1^{(p_1)}$ as a generator of
$\widetilde{\L}(1)\otimes\Z_{(p)}$ in the dimension $-2$
(i.e. in $\O_U^{-2}$).
Indeed, let $p_2$ be any prime. Then $[[x]_{p_2}]_{p_1}=[[x]_{p_1}]_{p_2}$.
Hence,
\begin{multline}
\label{pp}
  p_1p_2x+p_1\sum_n\a_n^{(p_2)}x^{n+1}+\sum_n\a_n^{(p_1)}
  \left(p_2x+\sum_m\a_m^{(p_2)}x^{m+1}\right)^{n+1}\\
  =p_2p_1x+p_2\sum_n\a_n^{(p_1)}x^{n+1}+\sum_n\a_n^{(p_2)}
  \left(p_1x+\sum_m\a_m^{(p_1)}x^{m+1}\right)^{n+1}.
\end{multline}
Taking the coefficient of $x^2$ in both sides of the above identity, we get
$p_1\a_1^{(p_2)}+p_2^2\a_1^{(p_1)}=p_2\a_1^{(p_1)}+p_1^2\a_1^{(p_2)}$.
Hence, $(p_1-p_1^2)\a_1^{(p_2)}=(p_2-p_2^2)\a_1^{(p_1)}$. Since $p_1$ is a
generator of $(\Z/p)^*$, the element $p_1-p_1^2$ is invertible in
$\Z_{(p)}$. So, $\a_1^{(p_2)}=\l\a_1^{(p_1)}$ with
$\l\in\Z_{(p)}\subset\O_U\otimes\Z_{(p)}$. Thus, for any prime
$p_2\ne p_1$ the coefficient $\a_1^{(p_2)}$ is a multiple of
$\a_1^{(p_1)}$, and that is why it can be excluded from the set of
generators for $\widetilde{\L}(1)\otimes\Z_{(p)}$.

Now, consider the coefficient system $\a_1,\ldots,\a_k,\ldots$ introduced
in the theorem. (That is, $\a_i$ is the coefficient of
$x^{i+1}$ in the series $[x]_{p_1}$ if $i$ is not divisible by $p-1$, and is
the coefficient of $x^{i+1}$ in the series $[x]_p$ if $i$ is divisible by
$p-1$.) By induction, we may suppose that this coefficient system is
a set of generators for $\widetilde{\L}(1)\otimes\Z_{(p)}$ in all
dimensions up to $-2(n-1)$. Hence, for any $q$ and $k\le n-1$ one has
\begin{equation}
\label{new}
  \a_k^{(q)}=\l_1^{(q)}\a_1+\ldots+\l_k^{(q)}\a_k
\end{equation}
with $\l_i^{(q)}\in\O_U\otimes\Z_{(p)}$. We are going to prove
that $\a_n^{(q)}$ can be also decomposed in such a way. It follows from
the above argument that we can consider only prime $q$.

First, suppose that $n$ is not divisible by $p-1$. Hence,
$\a_n=\a_n^{(p_1)}$, where $p_1$ is a generator of $(\Z/p)^*$.
Let $p_2$ be any prime. Taking the coefficient of $x^{n+1}$ in both
sides of~(\ref{pp}), we obtain
\begin{multline*}
  p_1\a_n^{(p_2)}+p_2^{n+1}\a_n^{(p_1)}+\m_1\a_1+\ldots+\m_{n-1}\a_{n-1}\\
  =p_2\a_n^{(p_1)}+p_1^{n+1}\a_n^{(p_2)}+\n_1\a_1+\ldots+\n_{n-1}\a_{n-1}
\end{multline*}
Here we expressed coefficients $\a_k^{(p_1)}$, $\a_k^{(p_2)}$, $k<n$, as a
linear combinations of generators $\a_1,\ldots,\a_{n-1}$, i.e.
$\m_i,\n_i\in\O_U\otimes\Z_{(p)}$. Therefore,
\begin{equation}
\label{p2}
  p_1(1-p_1^n)\a_n^{(p_2)}=(p_2-p_2^{n+1})\a_n^{(p_1)}+
  (\n_1-\m_1)\a_1+\ldots+(\n_{n-1}-\m_{n-1})\a_{n-1}.
\end{equation}
Since $p_1$ is a generator of $(\Z/p)^*$ and $n$ is not divisible by
$p-1$, we deduce that $p_1(1-p_1^n)$ is invertible in $\Z_{(p)}$. Thus,
it follows from~(\ref{p2}) that $\a_n^{(p_2)}$ is a linear
combination of $\a_1,\ldots,\a_{n-1}$ and $\a_n=\a_n^{(p_1)}$ with
coefficients from $\O_U\otimes\Z_{(p)}$.

Now, suppose that $n$ is divisible by $p-1$, i.e. $\a_n=\a_n^{(p)}$.
Before we proceed further, let us make some preliminary remarks.
It is well known (Milnor, Novikov), that the
complex cobordism coefficient ring $\O_U$ is a polynomial ring:
$\O_U=\Z[a_1,a_2,\ldots,a_n,\ldots]$, $a_n\in\O_U^{-2n}$. The ring
$\O_U$ is the coefficient ring of the (universal) formal
group law of geometric cobordisms (cf.~\cite{Q},~\cite{BuN}). This formal
group law has a logarithm series with coefficients in $\O_U\otimes\Q$, namely
$g(u)=u+\sum_n\frac{\C P^n}{n+1}u^{n+1}$, cf.~\cite{No1}.
Hence, the coefficient ring of the logarithm is
$\O_U(\Z):=\Z[b_1,b_2,\ldots,b_n,\ldots]$, where $b_n=\frac{\C P^n}{n+1}$.
It is well known that this ring is the maximal subring of $\O_U\otimes\Q$
on which all cohomological characteristic numbers take integer values.
One can choose generator sets $\{a_i^*\}$, $\{b_i^*\}$ for the rings
$\O_U$, $\O_U(\Z)$ such that the inclusion
$\iota_0:\O_U\to\O_U(\Z)$ is as follows:
$$
  \iota_0(a_i^*)=\left\{
  \begin{aligned}
    p\cdot b_i^*, & \quad\text{if $i=p^k-1$ for some $k>0$,}\\
    b_i^* & \quad\text{otherwise.}
  \end{aligned}\right.
$$
Let $B^+$ be the set of elements of degree $\ne0$ in the ring $B:=\O_U(\Z)$ .
Then $(B^+)^2$ consists of elements in $\O_U(\Z)$ that are decomposable
into the product of two non-trivial factors. The map
$\iota_0:\O_U\to\O_U(\Z)$ takes the coefficients $\a_n^{(p)}$ of the series
$[x]_p$ to the element of the form $(p-p^{n+1})b_n+((B^+)^2)$
(cf.~\cite[p.~22]{B2}). Therefore, the coefficients $\a_{p^k-1}^{(p)}$ can
be taken as multiplicative generators of $\O_U\otimes\Z_{(p)}$ in dimensions
$p^k-1$. In other dimensions $l\ne p^k-1$ we have
$\a_l^{(p)}\in p\O_U$, i.e. $\a_l^{(p)}$ is divisible by $p$ in $\O_U$.

Now, let us return to the proof of theorem~\ref{form}.
Let us rewrite the identity~(\ref{pp}) substituting $p$ for $p_1$:
\begin{multline*}
  pp_2x+p\sum_m\a_m^{(p_2)}x^{m+1}+\sum_m\a_m^{(p)}
  \left(p_2x+\a_1^{(p_2)}x^2+\a_2^{(p_2)}x^3+\ldots\right)^{m+1}\\
  =p_2px+p_2\sum_m\a_m^{(p)}x^{m+1}+\sum_m\a_m^{(p_2)}
  \left(px+\a_1^{(p)}x^2+\a_2^{(p)}x^3+\ldots\right)^{m+1}.
\end{multline*}
Taking the coefficient of $x^{n+1}$ in both sides, we get
\begin{multline*}
  p\a_n^{(p_2)}+p_2^{n+1}\a_n^{(p)}+\left< \sum_{m<n}\a_m^{(p)}
  (p_2x+\a_1^{(p_2)}x^2+\a_2^{(p_2)}x^3+\ldots)^{m+1} \right>_{n+1}\\
  =p_2\a_n^{(p)}+p^{n+1}\a_n^{(p_2)}+\left< \sum_{m<n}\a_m^{(p_2)}
  (px+\a_1^{(p)}x^2+\a_2^{(p)}x^3+\ldots)^{m+1} \right>_{n+1},
\end{multline*}
where $\langle\cdot\rangle_{n+1}$ denotes the coefficient of $x^{n+1}$.
Let us write again the coefficients $\a_m^{(p_2)}$ for
$m<n$ as linear combinations of generators $\a_1,\ldots,\a_m$. Since
$\a_m^{(p)}\in p\O_U$ for $m\ne p^k-1$, the last identity can be rewritten as
\begin{multline}
\label{np2}
  p(1-p^n)\a_n^{(p_2)}=p_2(1-p_2^n)\a_n^{(p)}
  +p(\m_1\a_1+\ldots+\m_{n-1}\a_{n-1})\\
  -\left< \sum_{k:p^k-1<n}\a_{p^k-1}^{(p)}
  \left( p_2x+\a_1^{(p_2)}x^2+\a_2^{(p_2)}x^3+\ldots \right)^{p^k}
  \right>_{n+1}\\
  +\left< \sum_{m<n}\a_m^{(p_2)}
  \left( \a_{p-1}^{(p)}x^p+\a_{p^2-1}^{(p)}x^{p^2}+\ldots+
  \a_{p^k-1}^{(p)}x^{p^k}+\ldots \right)^{m+1} \right>_{n+1}
\end{multline}
for some $\m_i\in\O_U\otimes\Z_{(p)}$.
The last two summands in the above formula can be rewritten as
$\a_{p-1}^{(p)}\n_1+\a_{p^2-1}^{(p)}\n_2+\ldots+\a_{p^k-1}^{(p)}\n_k$, where
$k=[\log_p(n+1)]$, $\n_i\in\O_U$. The
coefficients $\a_{p^i-1}^{(p)}$ are multiplicative generators of
$\O_U\otimes\Z_{(p)}$ in the dimensions $p^i-1$. Since $\O_U\otimes\Z_{(p)}$
is a polynomial ring, one has
$\n_i\in p\O_U\otimes\Z_{(p)}$, i.e. $\n_i$ is divisible by
$p$ in $\O_U\otimes\Z_{(p)}$. Let $\n_i=p\kappa_i$ with
$\kappa_i\in\O_U\otimes\Z_{(p)}$. Then~(\ref{np2}) gives
\begin{multline*}
  p(1-p^n)\a_n^{(p_2)}=p_2(1-p_2^n)\a_n^{(p)}+
  p(\m_1\a_1+\ldots+\m_{n-1}\a_{n-1})\\
  +p(\a_{p-1}^{(p)}\kappa_1+\a_{p^2-1}^{(p)}\kappa_2
  +\ldots+\a_{p^k-1}^{(p)}\kappa_k),
\end{multline*}
where $k=[\log_p(n+1)]$, $\m_i,\kappa_i\in\O_U\otimes\Z_{(p)}$. Since $n$ is
divisible by $p-1$, it follows that $1-p_2^n$ is divisible by $p$
(for $p_2\ne p$). Hence, the whole above identity is divisible by $p$.
Dividing it by $p$ and mentioning that $1-p^n$ is invertible in
$\O_U\otimes\Z_{(p)}$, we obtain that $\a_n^{(p_2)}$ is decomposable as
$$
  \a_n^{(p_2)}=\frac{p_2(1-p_2^n)}{p(1-p^n)}\a_n^{(p)}+\l_1\a_1+\ldots+
  \l_{n-1}\a_{n-1}
$$
with $\l_i\in\O_U\otimes\Z_{(p)}$.
Thus, setting
$$
  \l_n=\frac{p_2(1-p_2^n)}{p(1-p^n)}\:\in\:\O_U\otimes{\mathbb Z}_{(p)},
$$
we get a decomposition of type~(\ref{new}) for
$\a_n^{(p_2)}$ (note that $\a_n=\a_n^{(p)}$), which completes the
proof of theorem~\ref{form}.
\end{proof}

\begin{corollary}
\label{withp}
  Let $p_1$ be a prime generator of the cyclic group $(\Z/p)^*$.
  There is the following set of generators for the
  $\O_U\otimes\Z_{(p)}$-module $\widetilde{\L}_p(1)\otimes\Z_{(p)}$:
  $$
  \a_n=\left\{
  \begin{aligned}
    p, & \quad\text{if $n=0$,}\\
    \a_n^{(p_1)}, & \quad\text{if $n$ is not divisible by $p-1$,}\\
    \a^{(p)}_{p^k-1}, & \quad\text{if $n=p^k-1$, $k=1,2,\ldots$.}
  \end{aligned}\right.
  $$
  The $\O_U\otimes\Z/p$-module $\widetilde{\L}(1)\otimes\Z/p$
  has the following generator set:
  $$
  \a_n=\left\{
  \begin{aligned}
    \a_n^{(p_1)}, & \quad\text{if $n$ is not divisible by $p-1$,}\\
    \a^{(p)}_{p^k-1}, & \quad\text{if $n=p^k-1$, $k=1,2,\ldots$.}
  \end{aligned}\right.
  $$
\end{corollary}
\begin{remark}
  In both cases there no generators in the dimensions $n$ divisible by $p-1$
  other than $p^k-1$.
\end{remark}
\begin{proof}
Consider the set of generators for
$\widetilde{\L}(1)\otimes\Z_{(p)}$ constructed in theorem~\ref{form}.
If $n$ is divisible by $p-1$ and $n\ne p^k-1$, the
elements $\a_n$ are divisible by $p$, i.e. lie in $p\O_U$. All other
$\a_n$ do not belong to $p\O_U$.
\end{proof}

\section{Cohomological description of the set of cobordism classes of
manifolds with a simple action of $\Z/p$ and some corollaries}

In this section we use the previously obtained description of the
$\O_U\otimes\Z_{(p)}$-module $\widetilde{\L}_p(1)\otimes\Z_{(p)}$ to prove
the result analogous to the well-known Stong-Hattori theorem~\cite{CF2}.
Namely, we will describe the set of cobordism classes of manifolds
with a simple $\Z/p$-action in terms of the characteristic numbers.

As it was shown in~\cite{BuN}, the homomorphism
$\Phi:F_1\to\widetilde{\L}_p(1)\otimes\Z_{(p)}$ (see~(\ref{realiz}) and
following discussion) can be extended to a homomorphism
$\g_p:F_0\to\O_U(\Z)\otimes\Z_{(p)}$ such that
$$
  \g_p\bigl(\a_{2k-1}(x_1,x_2,\ldots,x_k)\bigr)=\left< \Bigl(
  \prod_{j=1}^k\frac u{[u]_{x_j}} \Bigr) \frac{pu}{[u]_p} \right>_k,
$$
where $\a_{2k-1}(x_1,x_2,\ldots,x_k)\in F_0$ is the Conner--Floyd invariant
(see~(\ref{rel2})). In particular,
$$
  \g_p\bigl(\a_{2k-1}(1,\ldots,1)\bigr)=
  \Bigl\langle \frac{pu}{[u]_p} \Bigr\rangle_k.
$$

Hence, for any simple action of $\Z/p$ on
$M^{2n}$ the $\mod p$ cobordism class of $M^{2n}$ can be expressed in terms
of the cobordism classes $\l_j\in\O_U$ of fixed submanifolds
and the weights $x_k^{(j)}\in(\Z/p)^*$ in the corresponding (trivial)
normal bundles as follows:
\begin{equation}
\label{class}
  [M^{2n}]\equiv\sum_j\l_j\g_p(x_1^{(j)},\ldots,x_{m_j}^{(j)}) \mod p\O_U.
\end{equation}
Now, the following question arises: which elements of the form
$$
  \sum_j\l_j\g_p(x_1^{(j)},\ldots,x_{m_j}^{(j)})\:\in\:\O_U({\mathbb
  Z})\otimes{\mathbb Z}_{(p)}
$$
are cobordism classes of manifolds with a simple
$\Z/p$-action? This question was firstly posed in~\cite{BuN} and is
analogous to the Milnor--Hirzebruch problem of describing the set of
elements in $\O_U({\Z})$ that are cobordism classes of (stably
complex) manifolds. While the Milnor--Hirzebruch problem is solved by the
Stong--Hattori theorem, the answer to the above question is given in our
theorem~\ref{charnum}. We will need the following definition.

\begin{definition}
  Let $\o=\sum_{i=1}^lk_i\cdot(i)$,
  $i,k_i\in\Z$, $i>0$, $k_i\ge0$, be a partition
  of $n=\|\o\|=\sum_ik_i\cdot i$ (i.e. $n$ is decomposed into the sum of
  positive integers, and the number $i$ enters this sum $k_i$ times). We say
  that the partition $\o$ is
  {\it divisible by $p-1$}, if all $i$ such that $k_i\ne 0$ are divisible by
  $p-1$ (i.e. all the summands are divisible by $p-1$; obviously, such
  partitions exist only for those $n$ divisible by $p-1$). We say that the
  partition $\o$ is {\it non $p$-adic}, if for any $j>0$ one has
  $k_{p^j-1}=0$ (i.e. there no summands of the form $p^j-1$).
\end{definition}

For each partition $\o=\sum_{i=1}^lk_i(i)$ let us put
$\|\o\|=\sum_ik_i\cdot i$ and $|\o|=\sum_ik_i$
(i.e. $\o$ is a partition of $\|\o\|$ with number of summands equals $|\o|$).
A~partition $\o$ defines a
characteristic class $s_\o$ as follows. Let us consider the smallest
symmetric polynomial in $x_1,\ldots,x_{|\o|}$ containing the monomial
$$
  (x_1\cdots x_{k_1})(x_{k_1+1}^2\cdots x_{k_1+k_2}^2)\cdots
  (x_{|\o|-k_l+1}^l\cdots x_{|\o|}^l).
$$
This polynomial defines a characteristic class in a usual way
(cf.~\cite{St}); in order to express it in terms of the Chern characteristic
classes one should write it as a polynomial in the elementary symmetric
functions $\sigma_i$, and then substitute $c_i$ for $\sigma_i$. Given an
$2n$-dimensional stably complex manifold $M^{2n}$, one can define the
cohomological characteristic numbers
$s_\o(M^{2n}):=s_\o(TM^{2n})[M^{2n}]\in\Z$ for all partitions $\o$ such
that $\|\o\|=n$. The $K$-theory characteristic numbers $s_\o(M^{2n})\in\Z$
(cf.~\cite{CF2}) are defined for all partitions $\o$ such that
$\|\o\|\le n$; they coincide with the cohomological numbers for
$\|\o\|=n$, while for $\o=0$ the corresponding $K$-theory characteristic
number is the Todd genus.

\begin{theorem}
\label{charnum}
  An element
  $\s\in\O_U(\Z)^{-2n}\otimes\Z_{(p)}$ belongs to the
  $\O_U\otimes\Z_{(p)}$-module $\widetilde{\L}_p(1)\otimes\Z_{(p)}$ and
  therefore, is the cobordism class of a manifold with a simple
  $\Z/p$-action, if and only if all its $K$-theory characteristic numbers
  $s_{\o}(\s)$, $\o=\sum_ik_i\cdot(i)$,
  $\|\o\|=\sum_ik_i\cdot i\le n$, lie in $\Z_{(p)}$, and for all
  partitions $\o$ divisible by $p-1$ the cohomological characteristic
  numbers $s_{\o}(\s)$, $\|\o\|=n$, are zero modulo $p$.
\end{theorem}

\begin{proof} (a) Necessity.\\
Let $\s\in\widetilde{\L}_p(1)\otimes\Z_{(p)}$. Note that the set
of generators for the $\O_U\otimes\Z_{(p)}$-module
$\widetilde{\L}_p(1)\otimes\Z_{(p)}$ described in corollary~\ref{withp}
has the following property: each of its elements
$\a_i\in\widetilde{\L}_p(1)\otimes\Z_{(p)}$ is a
multiplicative generator of $\O_U\otimes\Z_{(p)}$ in dimension $-2i$.
However, this set of generators for $\widetilde{\L}_p(1)\otimes\Z_{(p)}$
has no elements in dimensions $-2i$ such that $i$ is divisible by $p-1$ and
$i\ne p^k-1$. So, we add any generators $\a_i$ in these missing
dimensions to get the whole set of multiplicative generators for
$\O_U\otimes\Z_{(p)}$. Now we have
$\O_U\otimes\Z_{(p)}=\Z_{(p)}[\a_1,\a_2,\ldots]$.

Since $\s\in\widetilde{\L}_p(1)\otimes\Z_{(p)}
\subset\O_U\otimes\Z_{(p)}$, it follows from the Stong--Hattori theorem that
all the $K$-characteristic numbers $s_{\o}(\s)$, $\|\o\|\le n$, lie in
$\Z_{(p)}$.

If $n$ is not divisible by $p-1$, then there are no partitions $\o$
divisible by $p-1$.

Now, let $n=m(p-1)$. One can write $\s$ as a homogeneous polynomial of
degree $-2m(p-1)$ in $\a_i$:
\begin{equation}
\label{sigma}
  \s=\sum_{\|\o\|=m(p-1)}r_{\o}\a_{\o}=r_{m(p-1)}\a_{m(p-1)}+\ldots,
\end{equation}
where $\a_{\o}=\a_1^{k_1}\cdot\a_2^{k_2}\cdots\a_l^{k_l}$ for
$\o=\sum_ik_i\cdot(i)$. It follows from the description of
$\widetilde{\L}_p(1)\otimes\Z_{(p)}$ given in corollary~\ref{withp} that
$\s\in\widetilde{\L}_p(1)\otimes\Z_{(p)}$ if and only if the coefficients
$r_{\o}$ in the decomposition~(\ref{sigma}) are zero modulo $p$ for all non
$p$-adic and divisible by $p-1$ partitions $\o$.

Consider the Chern--Dold character $\mathop{\rm ch}_U:U^{*}(\:\cdot\:)\to
H^{*}(\:\cdot\:;\O_U\otimes\Q)$ in cobordisms~\cite{Bu}:
$$
  \mathop{\rm ch}\nolimits_U(v)=t+\sum_{i\ge 1}\beta_it^{i+1}.
$$
Here $v=c_1^{U}(\zeta)\in U^2(\C P^\infty)$ is the
first cobordism Chern class of the universal line bundle,
$t=c_1^{H}(\zeta)\in H^2(\C P^\infty)$ is the same Chern class in
cohomologies, and the coefficients $\beta_i$ are from $\O_U(\Z)$. Then for
any $\s\in\O_U^{-2n}$ holds
\begin{equation}
\label{cn}
  \s=\sum_{\|\o\|=n}s_{\o}(\s)\beta_{\o},
\end{equation}
where $\beta_{\o}=\beta_1^{k_1}\cdot\beta_2^{k_2}\cdots\beta_l^{k_l}$ for
$\o=\sum_ik_i\cdot(i)$. The coefficient ring
${\mathbb Z}[\beta_1,\beta_2,\ldots]$ of
the Chern--Dold character coincides with $B=\O_U(\Z)$ (cf.~\cite{Bu}).
Hence,
$$
\a_i=\left\{
\begin{aligned}
  e_i\cdot\beta_i+((B^+)^2) &\quad \text{if $i\ne p^k-1$},\\
  pe_i\cdot\beta_i+p((B^+)^2) &\quad \text{if $i=p^k-1$}
\end{aligned}\right.
$$
with invertible
$e_i\in\Z_{(p)}$. Now, let us write $\s$ as a homogeneous polynomial
in $\beta_i$. Since all $\beta_i\in\O_U(\Z)$
have integer cohomological characteristic
numbers, to prove the necessity of the theorem it suffices to show that
the coefficient of $\beta_{\o}$ in the decomposition of $\s$ is zero
modulo $p$ if the partition $\o=\sum_ik_i\cdot(i)$ is divisible by $p-1$.
This coefficient is the homological characteristic number
$s_{\o}(\s)$ (see~(\ref{cn})), which can be decomposed as follows
(see~(\ref{sigma})):
\begin{equation}
\label{coef}
  s_{\o}(\s)=\sum_{\o^{\prime}:\o^{\prime}\supset\o}
  r_{\o^{\prime}}s_{\o}(\a_{\o^{\prime}}),
\end{equation}
where $\o^{\prime}\supset\o$ means that $\o$ refines $\o'$.
This coefficient is divisible by $p$. Indeed, if the partition
$\o'=\sum_ik^{\prime}_i\cdot(i)$ is divisible by $p-1$ and
non $p$-adic, then $r_{\o^{\prime}}$ is zero modulo $p$, since
$\s\in\widetilde{\L}_p(1)\otimes\Z_{(p)}$ (see above). If there are
some summands of the form $p^k-1$ in the partition $\o'$, then
$\a_{\o^{\prime}}\in p\O_U(\Z)\otimes\Z_{(p)}$, i.e.
$s_{\o}(\a_{\o^{\prime}})$ is divisible by $p$. Anyway, the whole sum
in~(\ref{coef}) is divisible by $p$. The necessity of the theorem is proved.

\medskip

(b) Sufficiency.\\
Since all the $K$-characteristic numbers of $\s$ are in $\Z_{(p)}$, it
follows from the Stong--Hattori theorem~\cite{CF2} that
$\s\in\O_U\otimes\Z_{(p)}$. Besides, suppose that the
characteristic numbers $s_{\o}(\s)$ are zero modulo $p$ for all divisible by
$p-1$ partitions $\o=\sum_ik_i\cdot(i)$, $\|\o\|=n$ .

Consider again the constructed above generator set $\a_1,\a_2,\ldots$ for
$\O_U\otimes\Z_{(p)}$. In order to prove that
$\s\in\widetilde{\L}_p(1)\otimes\Z_{(p)}$ one needs to show that for every
divisible by $p-1$ and non $p$-adic partition $\o=\sum_ik_i\cdot(i)$ the
coefficient $r_{\o}$ in decomposition~(\ref{sigma}) is zero modulo $p$.
Let $\o$ be such a partition. We can rewrite identity~(\ref{coef}) as
follows:
\begin{equation}
\label{coef2}
  s_{\o}(\s)=r_{\o}s_{\o}(\a_{\o})+
  \sum_{\o^{\prime}\supset\o,\o^{\prime}\ne\o}
  r_{\o^{\prime}}s_{\o}(\a_{\o^{\prime}}).
\end{equation}

One can assume by induction that if a
partition $\o'$ such that $\o^{\prime}\supset\o$, $\o^{\prime}\ne\o$,
$\|\o^{\prime}\|=m(p-1)$, is non $p$-adic, then
the coefficient $r_{\o^{\prime}}$ is
divisible by $p$. If the partition $\o'=\sum_ik_i^{\prime}\cdot(i)$ is not
non $p$-adic (i.e. there some summands of the form $p^k-1$), then
$s_{\o}(\a_{\o^{\prime}})$ is divisible by $p$. Anyway, the second
summand in the right hand side of~(\ref{coef2}) is zero modulo $p$.
The left hand side of~(\ref{coef2}) is zero modulo $p$ by assumption. Since
$\o$ is non $p$-adic, we have
$\a_{\o}=e\cdot\beta_{\o}+\ldots$ with invertible $e\in\Z_{(p)}$. So,
$s_{\o}(\a_{\o})$ is not divisible by $p$. Thus, it follows
from~(\ref{coef2}) that $r_{\o}$ is zero modulo $p$.
\end{proof}

\begin{corollary}
\label{cobord}
  An element $\s\in\O_U$ is the cobordism class of a manifold with
  $\Z/p$-action whose fixed point set has the trivial normal bundle if and
  only if the cohomological
  characteristic numbers $s_{\o}(\s)$, $\|\o\|=n$, are zero modulo $p$
  for all divisible by $p-1$ partitions $\o$.
\end{corollary}

\begin{corollary}
\label{smdim}
  Each cobordism class of dimension $n\le4p-6$ contains a manifold
  $M^n$ with a simple action of $\Z/p$.
\end{corollary}
This result was firstly proved in~\cite{Ka}.
In dimension $n=4p-4$ there exist a manifold
(e.g. $\C P^{2p-2}$) whose cobordism class does not contain a manifold
with a simple action of $\Z/p$.

In Conner and Floyd's book~\cite{CF1} it was shown by the methods not
involving the formal group theory, that a cobordism class $\sigma\in\O_U$
contains a manifold with a {\it strictly} simple action of $\Z/p$ (see
definition~\ref{simple}) if and only if {\it all} the characteristic numbers
$\s_{\o}(\s)$ are zero modulo $p$. More precisely, it was shown there
that the set of cobordism classes of manifolds with a strictly simple
$\Z/p$-action coincides with the $\O_U$-module spanned by the set
$Y^0=p,Y^1,Y^2,\ldots$, where $Y^i\in\O_U^{p^i-1}$ are the so-called ``Milnor
manifolds". These manifolds $Y^i$ are uniquely determined by the following
conditions: $s_{(p^i-1)}(Y^i)=p$, and $s_{\o}(Y^i)$ is divisible by $p$ for
any $\o$. For our purposes we may consider
$\O_U\otimes\Z_{(p)}$-modules instead of $\O_U$-modules. Hence, one could
take the elements $\a_{p^i-1}^{(p)}$ from corollary~\ref{withp} as
representatives of the cobordism classes of $Y^i$. Now, we see that the
$\O_U\otimes\Z_{(p)}$-module $\O_U[p,Y^1,Y^2,\ldots]\otimes\Z_{(p)}$ studied
by Conner and Floyd is included into our $\O_U\otimes\Z_{(p)}$-module
$\widetilde{\L}_p(1)\otimes\Z_{(p)}$, and the set of generators for the
former module is a subset of the generator set for the latter one.

Finally, we note that if a certain cobordism class $\sigma\in\O_U$
contains a representative $M$ with a strictly simple action of $\Z/p$, then
it is not necessarily true that {\it any} simple action of $\Z/p$ on $M$ is
strictly simple. Indeed, let us consider two simple actions, first on $M_1=\C
P^{p-1}$ with generator $\rho\in\Z/p$ acting as
$\rho(z_1:\ldots:z_p)=(z_1:\rho z_2:\ldots:\rho^{p-1}z_p)$ (this simple
action with $p$ fixed points is strictly simple as well), and second on
$M_2=\C P^{1}$, $\rho(z_1:z_2)=(z_1:\rho z_2)$ (this simple action with $2$
fixed points is not strictly simple). Then one has two simple $\Z/p$-actions
on $M=M_1\times M_2$: $\rho(a,b)=(\rho a,b)$ and $\rho(a,b)=(a,\rho
b)$, $a\in\C P^{p-1}$, $b\in\C P^{1}$. The first one is strictly simple,
while the second one is not.

The author is grateful to Prof. V.\,M.~Buchstaber for useful
recommendations, stimulating discussions and attention to the research.


\begin{thebibliography}{CF2}

\bibitem[B1]{Bu}
V.\,M. Bukhshtaber,
{\it The Chern--Dold character in cobordisms.~I,}
Mat. Sbornik {\bf 83} (1970), no.~4, 575--595;
English transl. in Math. USSR, Sbornik {\bf 12} (1970), 573--594.

\bibitem[B2]{B2}
V.\,M. Bukhshtaber,
{\it Characteristic classes in cobordisms and topological applications of
one-valued and two-valued formal group law theories,}
Sovremennye Problemy Matematiki {\bf 10} (1978), 5--178 (in Russian).

\bibitem[BN]{BuN}
V.\,M. Bukhshtaber and S.\,P.Novikov,
{\it Formal groups, power systems and Adams operators,}
Mat. Sbornik {\bf 84} (1971), no.~1, 81--118;
English transl. in Math. USSR, Sbornik {\bf 13} (1971), 80--116.

\bibitem[CF1]{CF1}
P.\,E. Conner and E.\,E. Floyd,
{\it Differentiable periodic maps,}
Springer-Verlag, Berlin, 1964.

\bibitem[CF2]{CF2}
P.\,E. Conner and E.\,E. Floyd,
{\it The relation of cobordism to K-theories,}
Lecture Notes in Math., Springer-Verlag, Berlin-Heidelberg-New York, 1966.

\bibitem[K]{Ka}
G.\,G. Kasparov,
{\it Invariants of classical lens manifolds in cobordism theory,}
Izvestiya Akad. Nauk SSSR Ser. Mat. {\bf 33} (1969), no.~4, 735--747;
English transl. in Math. USSR, Izvestiya {\bf 3} (1969), 695--705.

\bibitem[N1]{No1}
S.\,P. Novikov,
{\it The methods of algebraic topology from the viewpoint of cobordism
theory,}
Izvestiya Akad. Nauk SSSR Ser. Mat. {\bf 31} (1967), no.~4, 855--951;
English transl. in Math. USSR, Izvestiya {\bf 1} (1967), 827--913.

\bibitem[N2]{No2}
S.\,P. Novikov,
{\it Adams operators and fixed points,}
Izvestiya Akad. Nauk SSSR Ser. Mat. {\bf 32} (1968), no.~6, 1245--1263;
English transl. in Math. USSR, Izvestiya {\bf 2} (1968), 1193--1211.

\bibitem[Q]{Q}
D. Quillen,
{\it On the formal group laws of unoriented and complex cobordism theory,}
Bull. Amer. Math. Soc. {\bf 75} (1969), no.~6, 1293--1298.

\bibitem[S]{St}
R.\,E. Stong,
{\it Notes on cobordism theory,}
Princeton Univ. Press, Princeton, NJ, 1968.

\end{thebibliography}
\end{document}